\documentclass[11pt]{amsart}

\usepackage[margin=1in]{geometry}
\usepackage[utf8]{inputenc}
\usepackage[T1]{fontenc}
\usepackage{amsmath,amssymb,amsthm}
\usepackage{graphicx}
\usepackage{booktabs}
\usepackage{xcolor}
\usepackage{fancyvrb}
\usepackage{newunicodechar}
\usepackage{xurl}
\usepackage[colorlinks=true,linkcolor=blue!50!black,citecolor=blue!50!black,urlcolor=blue!50!black]{hyperref}

\newunicodechar{ℕ}{\ensuremath{\mathbb{N}}}
\newunicodechar{ℂ}{\ensuremath{\mathbb{C}}}
\newunicodechar{→}{\ensuremath{\to}}
\newunicodechar{↔}{\ensuremath{\leftrightarrow}}
\newunicodechar{∀}{\ensuremath{\forall}}
\newunicodechar{∃}{\ensuremath{\exists}}
\newunicodechar{∈}{\ensuremath{\in}}
\newunicodechar{∨}{\ensuremath{\lor}}
\newunicodechar{∧}{\ensuremath{\land}}
\newunicodechar{₁}{\ensuremath{{}_1}}
\newunicodechar{₂}{\ensuremath{{}_2}}
\fvset{fontsize=\small,xleftmargin=1em}

\newtheorem{theorem}{Theorem}[section]
\newtheorem{lemma}[theorem]{Lemma}
\newtheorem{proposition}[theorem]{Proposition}
\newtheorem{corollary}[theorem]{Corollary}

\theoremstyle{definition}

\newtheorem{example}[theorem]{Example}
\theoremstyle{remark}
\newtheorem{remark}[theorem]{Remark}

\newcommand{\F}{\mathbb{F}_2}
\newcommand{\rev}{\operatorname{rev}}
\newcommand{\tc}{\tilde c}

\newif\ifwithdoubling
\withdoublingtrue
\newcommand{\furthersections}{\ifwithdoubling Sections~\ref{sec:cf} and~\ref{sec:doubling}\else Section~\ref{sec:cf}\fi}

\begin{document}

\title[Nilpotent skew-symmetric tridiagonal Bohemian matrices]{Nilpotent skew-symmetric tridiagonal Bohemian matrices: a proof of Corless's conjecture}

\author{Marcel Blattner}
\address{Applied AI Research Lab, Lucerne University of Applied Sciences and Arts, Suurstoffi~1, 6343 Rotkreuz, Switzerland}
\email{marcel.blattner@hslu.ch}
\urladdr{https://orcid.org/0009-0001-7188-673X}


\begin{abstract}
Corless studied the skew-symmetric tridiagonal matrices $T_m(u)$ with superdiagonal $u \in \{1, i\}^{m-1}$. He showed that $T_m(u)$ can be nilpotent only if $m = 2^K - 1$ and that the $2^{K-1}$ vectors $u$ obtained from the empty vector by repeated application of $u \mapsto (u, 1, i, \rev u)$ and $u \mapsto (u, i, 1, \rev u)$ give nilpotent matrices, and he conjectured that there are no others. We prove this conjecture. Reducing the characteristic polynomial modulo $4$ gives a necessary condition for nilpotency that is linear over $\F$. For $m > 1$ it forces $u = (v, 1, i, \rev v)$ or $u = (v, i, 1, \rev v)$, and an identity implicit in Corless's proof then reduces the problem to $v$. Since $T_m(u)$ is diagonally similar to a finite section of the Feinberg--Zee hopping matrix, the theorem also determines the nilpotent tridiagonal sign matrices. We further show that the $(1,1)$ entry of the resolvent of a nilpotent $T_m(u)$ is $\lambda^{-1} \pm \lambda^{-3} \pm \lambda^{-7} \pm \dots \pm \lambda^{-(2^K - 1)}$, with signs that determine $u$. The proof of the main theorem has been formalised in Lean~4.
\end{abstract}

\maketitle

\section{Introduction}

A matrix family is called \emph{Bohemian} if its entries are drawn from a fixed finite set, its \emph{population} \cite{Cor21}. Corless \cite{Cor21} studied the skew-symmetric tridiagonal Bohemian matrices with population $\{1, i\}$:
\[
T_m(u) =
\begin{pmatrix}
0 & u_1 & & & \\
-u_1 & 0 & u_2 & & \\
 & -u_2 & 0 & \ddots & \\
 & & \ddots & \ddots & u_{m-1} \\
 & & & -u_{m-1} & 0
\end{pmatrix},
\qquad u = (u_1, \dots, u_{m-1}) \in \{1, i\}^{m-1}.
\]
Here skew-symmetric means $T^{\mathsf{T}} = -T$, without complex conjugation. Nilpotent matrices are of particular interest in this family, because their computed eigenvalues are far from the true ones: a floating-point eigenvalue routine returns values of size about $\mu^{1/m}$, where $\mu$ is the unit roundoff, and these values produce conspicuous artefacts in plots of the eigenvalues of the whole family \cite[Section~3]{Cor21}. A nilpotent $T_m(u)$ consists of a single Jordan block (Remark~\ref{rem:jordan}); Corless calls such matrices \emph{coreless}.

Corless proved that $T_m(u)$ can be nilpotent only if $m = 2^K - 1$ for some $K \ge 1$ \cite[Theorem~4.1]{Cor21}, and that the vectors obtained from the empty vector by repeated application of
\begin{equation}\label{eq:rule}
u \mapsto (u, 1, i, \rev u) \qquad\text{and}\qquad u \mapsto (u, i, 1, \rev u)
\end{equation}
give nilpotent matrices \cite[Theorem~4.2]{Cor21}; here $\rev u$ denotes $u$ reversed. This construction yields $2^{K-1}$ nilpotent matrices of dimension $2^K - 1$. He conjectured that there are no others \cite[Conjecture~5.1]{Cor21} and verified this by exhaustive computation for $m \le 31$.

\begin{theorem}\label{thm:main}
\interlinepenalty=10000 % keep the statement on one page
Let $m \ge 1$ and $u \in \{1, i\}^{m-1}$. Then $T_m(u)$ is nilpotent if and only if $m = 2^K - 1$ for some $K \ge 1$ and $u$ is obtained from the empty vector by repeated application of \eqref{eq:rule}. In particular, there are exactly $2^{K-1}$ nilpotent matrices $T_m(u)$ of dimension $m = 2^K - 1$, and none of any other dimension.
\end{theorem}

\begin{figure}[t]
\centering
\includegraphics[width=0.95\textwidth]{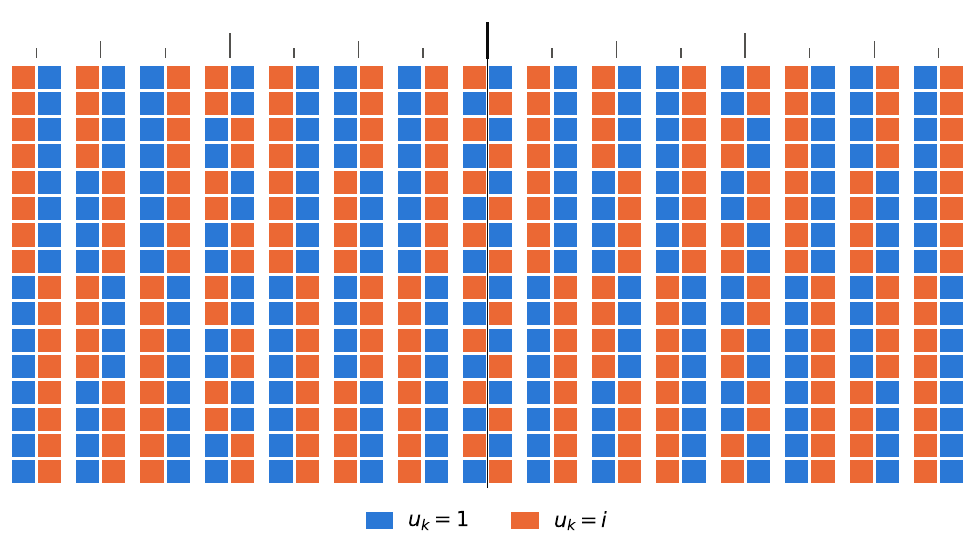}
\caption{The $16$ nilpotent vectors $u$ for $m = 31$, one per row; blue cells are $1$, orange cells are $i$. Consecutive entries are grouped in pairs (see Corollary~\ref{cor:pairs}), and the ticks mark the middle pair of each level of the recursion.}
\label{fig:family}
\end{figure}

For $m = 3$ the nilpotent vectors are $(1, i)$ and $(i, 1)$; Figure~\ref{fig:family} shows the $16$ vectors for $m = 31$.

\subsection{Outline of the proof}

Put $s_k = u_k^2 = \pm 1$. The characteristic polynomial of $T_m(u)$ is $\lambda^m D_s(\lambda^{-2})$, where $D_s$ is the matching polynomial of the path with $m$ vertices and edge weights $s_1, \dots, s_{m-1}$; see \cite[Propositions~2.1 and~2.9]{Cor21}. Thus $T_m(u)$ is nilpotent if and only if $D_s = 1$, which is a system of polynomial equations of degrees $1, 2, \dots, \lfloor m/2 \rfloor$ in the signs. Call a sign sequence \emph{folded} if it has the form $s = (\sigma, \varepsilon, -\varepsilon, \rev\sigma)$ with $\varepsilon = \pm 1$; the rule \eqref{eq:rule} produces such sequences. Folded sequences satisfy the identity $D_s = D_\sigma^2$ (Lemma~\ref{lem:fold}), which is implicit in the proof of \cite[Theorem~4.2]{Cor21}. Hence a folded sequence satisfies $D_s = 1$ if and only if $D_\sigma = 1$, and the theorem follows by induction once we know that every nonempty $s$ with $D_s = 1$ is folded.

To prove this we study the equation $D_s = 1$ modulo $2$ and modulo $4$. Modulo $2$ the signs disappear, and one obtains the restriction on the dimension. Modulo $4$ they enter linearly: if $s_e = 1 - 2t_e$ with $t_e \in \{0, 1\}$, then $\prod_e (1 - 2t_e) \equiv 1 - 2 \sum_e t_e \pmod 4$. In other words, modulo $4$ the map $s \mapsto D_s$ coincides with its first-order Taylor expansion at $s = (1, \dots, 1)$, with the element $2$ of $\mathbb{Z}/4\mathbb{Z}$ in the role of an infinitesimal of square zero. The condition $D_s \equiv 1 \pmod 4$ is therefore a linear system over $\F$ in the unknowns $t_e$. By the reflection symmetry of the path, the edges $e$ and $m - e$ have the same coefficient polynomial, so the system involves only the sums $t_e + t_{m-e}$ in $\F$. Apart from this symmetry its coefficient polynomials are linearly independent (Lemma~\ref{lem:indep}), so the system determines these sums. Since every folded sequence is a solution, the sums are those of a folded sequence, and every solution is folded (Proposition~\ref{prop:mod4}). The moduli $8, 16, \dots$, for which the equations are no longer linear, are not needed.

\subsection{Tridiagonal sign matrices}\label{sec:sign}

The matrices $T_m(u)$ also occur, in a different normalisation, in the spectral theory of non-self-adjoint operators. For $a, c \in \mathbb{C}^{m-1}$ let $J_m(a, c)$ be the $m \times m$ tridiagonal matrix with zero diagonal, superdiagonal $a$ and subdiagonal $c$, so that $T_m(u) = J_m(u, -u)$. For $b \in \{\pm 1\}^{m-1}$ put $A_m^b = J_m\big((1, \dots, 1), b\big)$. These are the finite sections of the hopping matrix of Feinberg and Zee \cite{FZ99}, whose spectrum has been studied in \cite{CCL11, CCL13, CWD12, Hag15a}. Following Hagger \cite{Hag15b}, we call them \emph{tridiagonal sign matrices}; every tridiagonal matrix with zero diagonal and entries $\pm 1$ on the two off-diagonals is diagonally similar to one of them. So is $T_m(u)$, with $b_k = -u_k^2$ (Lemma~\ref{lem:charpoly}), and Theorem~\ref{thm:main} therefore classifies the nilpotent tridiagonal sign matrices as well. To state the result in this setting, let $F_K \subseteq \{\pm 1\}^{2^K - 2}$ be the set of sign sequences obtained from the empty sequence by $K - 1$ applications of $\sigma \mapsto (\sigma, \varepsilon, -\varepsilon, \rev\sigma)$ with $\varepsilon = \pm 1$; under $s_k = u_k^2$ these are the vectors of Theorem~\ref{thm:main}.

\begin{theorem}\label{thm:sign}
Let $m \ge 1$ and let $a, c \in \mathbb{C}^{m-1}$ satisfy $a_k c_k = \pm 1$ for all $k$. Then $J_m(a, c)$ is nilpotent if and only if $m = 2^K - 1$ for some $K \ge 1$ and $(a_1 c_1, \dots, a_{m-1} c_{m-1}) \in F_K$. In particular, a tridiagonal sign matrix $A_m^b$ is nilpotent if and only if $m = 2^K - 1$ for some $K \ge 1$ and $b \in F_K$.
\end{theorem}

Some nilpotent tridiagonal sign matrices are implicit in this literature. Chandler-Wilde, Chonchaiya and Lindner \cite{CCL11} introduced the sign sequence
\begin{equation}\label{eq:ctilde}
\tc_1 = 1, \qquad \tc_{2j} = \tc_{2j-1}\, \tc_j, \qquad \tc_{2j+1} = -\tc_{2j} \qquad (j \ge 1)
\end{equation}
(we use the notation of \cite{CWD12}) and showed that the characteristic polynomials of the finite sections $A_m^b$ with $b = (\tc_2, \dots, \tc_m)$ have all their coefficients in $\{-1, 0, 1\}$. Chandler-Wilde and Davies showed that for $m = 2^n - 1$ these polynomials are monomials \cite[Lemma~30]{CWD12}, so that the corresponding matrices $A_m^b$ are nilpotent. Together with the symmetries $b \mapsto -b$ \cite[Lemma~3.4]{CCL13} and $b \mapsto \rev b$, this accounts for four of the $2^{n-1}$ nilpotent tridiagonal sign matrices of dimension $2^n - 1 \ge 7$. \ifwithdoubling We give another proof in Remark~\ref{rem:ctilde}. \fi To the best of our knowledge, neither the restriction on the dimension nor the classification appears in this literature. The resemblance between the eigenvalue plots of the two families was pointed out in \cite{CLPS22}, where the ring at the centre of \cite[Figure~3]{CCL13} for dimension $15$ is attributed to rounding errors caused by nilpotent matrices.

\subsection{Further results}

In Section~\ref{sec:cf} we show that if $T_m(u)$ is nilpotent, then the terminating continued fraction with partial numerators $u_1^2 x, \dots, u_{m-1}^2 x$ is a lacunary polynomial $1 \pm x \pm x^3 \pm x^7 \pm \dots \pm x^{2^{K-1}-1}$, whose signs are the choices made in the recursion. Equivalently, the $(1,1)$ entry of the resolvent of $T_m(u)$ is $\lambda^{-1} \pm \lambda^{-3} \pm \lambda^{-7} \pm \dots \pm \lambda^{-(2^K-1)}$. \ifwithdoubling In Section~\ref{sec:doubling} we show that the nilpotent sign sequences are also generated by a doubling rule, of which the sequence \eqref{eq:ctilde} is a fixed point. This rule is a finite analogue of the maps used in \cite{CWD12} to show that the spectrum of the Feinberg--Zee hopping matrix is closed under taking square roots. \fi The proof of Theorem~\ref{thm:main} has been formalised in Lean~4 \cite{dMU21} with the Mathlib library \cite{mathlib}; this and an exhaustive search for $m = 63$ are described in Section~\ref{sec:lean}.

\section{Matching polynomials}\label{sec:matching}

Throughout, $s_k = u_k^2$, so $s_k = 1$ if $u_k = 1$ and $s_k = -1$ if $u_k = i$. Since $u_k$ is determined by $s_k$, we identify a vector $u$ with its \emph{sign sequence} $s = (s_1, \dots, s_{m-1}) \in \{\pm 1\}^{m-1}$. Corless's rule \eqref{eq:rule} becomes
\begin{equation}\label{eq:rule-s}
\sigma \mapsto (\sigma, \varepsilon, -\varepsilon, \rev \sigma), \qquad \varepsilon = \pm 1,
\end{equation}
so the vectors of Theorem~\ref{thm:main} correspond to the sequences in $F_K$, which have length $2^K - 2$. Since $\sigma$ and $\varepsilon$ can be read off $(\sigma, \varepsilon, -\varepsilon, \rev\sigma)$, we have $|F_K| = 2^{K-1}$, and $F_K$ is closed under $s \mapsto -s$, by induction on $K$.

Let $P_m$ be the path with the vertices $1, \dots, m$, in which the edge $e$ joins $e$ and $e + 1$ ($1 \le e \le m - 1$). For $s \in \{\pm 1\}^{m-1}$ let
\[
D_s(x) = \sum_{M} \Big(\prod_{e \in M} s_e\Big)\, x^{|M|} \in \mathbb{Z}[x],
\]
where $M$ runs through the matchings of $P_m$, that is, the sets of pairwise disjoint edges, including the empty set. For $s = (1, \dots, 1)$ we obtain the Fibonacci polynomials $f_m(x) = \sum_j \binom{m-j}{j} x^j$. More generally, for $1 \le a \le b + 1 \le m + 1$ let $D_{a,b}$ be the matching polynomial of the subpath with the vertices $a, \dots, b$ and the signs $s_a, \dots, s_{b-1}$. Thus $D_s = D_{1,m}$, and $D_{a,b} = 1$ if the subpath has at most one vertex. Distinguishing whether or not the last vertex is matched, we obtain
\begin{equation}\label{eq:rec}
D_{a,b} = D_{a,b-1} + s_{b-1}\, x\, D_{a,b-2} \qquad (b \ge a + 1).
\end{equation}

\begin{lemma}\label{lem:charpoly}
Let $a, c \in \mathbb{C}^{m-1}$ with $a_k c_k = \pm 1$ for all $k$.
\begin{enumerate}
\item[(i)] Put $s_k = -a_k c_k$. Then $\det(\lambda I - J_m(a, c)) = \lambda^m D_s(\lambda^{-2})$. In particular, $J_m(a, c)$ is nilpotent if and only if $D_s = 1$.
\item[(ii)] Put $b_k = a_k c_k$, and let $\Delta$ be the diagonal matrix with the entries $\delta_1 = 1$ and $\delta_{k+1} = \delta_k / a_k$. Then $J_m(a, c) = \Delta A_m^b \Delta^{-1}$.
\end{enumerate}
\end{lemma}

\begin{proof}
The characteristic polynomials $p_k$ of the leading $k \times k$ blocks of $J_m(a, c)$ satisfy $p_k = \lambda p_{k-1} - a_{k-1} c_{k-1}\, p_{k-2}$ with $p_0 = 1$ and $p_1 = \lambda$ (cf.\ \cite[Proposition~2.1]{Cor21}). By \eqref{eq:rec}, so do the polynomials $\lambda^k D_{1,k}(\lambda^{-2})$. This proves (i), and (ii) is verified entrywise.
\end{proof}

For $T_m(u) = J_m(u, -u)$ the lemma gives $s_k = u_k^2$ and $b_k = -u_k^2$. We call a vector $u$, and its sign sequence $s$, \emph{nilpotent} if $T_m(u)$ is nilpotent, that is, if $D_s = 1$. Theorem~\ref{thm:main} is equivalent to the following statement.

\begin{theorem}\label{thm:main-s}
Let $m \ge 1$ and $s \in \{\pm 1\}^{m-1}$. Then $D_s = 1$ if and only if $m = 2^K - 1$ and $s \in F_K$.
\end{theorem}

The rule \eqref{eq:rule-s} is explained by the following identity.

\begin{lemma}\label{lem:fold}
Let $\sigma \in \{\pm 1\}^{N-1}$ and $\varepsilon = \pm 1$, and put $s = (\sigma, \varepsilon, -\varepsilon, \rev\sigma)$. Then $D_s = D_\sigma^2$.
\end{lemma}

\begin{proof}
Deleting the middle vertex of the path of $s$ leaves two paths that are mirror images of each other, each with matching polynomial $D_\sigma$. Let $E$ be the matching polynomial of either of them without the vertex next to the middle one. Classifying matchings by how the middle vertex is covered (not at all, by the left edge, or by the right edge) gives
\[
D_s = D_\sigma^2 + \varepsilon x\, E D_\sigma - \varepsilon x\, D_\sigma E = D_\sigma^2. \qedhere
\]
\end{proof}

For the matrices $T_m(u)$ this identity is implicit in the proof of \cite[Theorem~4.2]{Cor21}: the first line of equation~(17) there comes from a Laplace expansion along the central row and does not use nilpotency.

Modulo $2$ all sign sequences give the same polynomial.

\begin{lemma}\label{lem:mod2}
Let $m \ge 1$ and $s \in \{\pm 1\}^{m-1}$.
\begin{enumerate}
\item[(i)] $D_s \equiv f_m \pmod 2$.
\item[(ii)] In $\F[y]$ we have $f_k(y + y^2) = (1 + y)^{k+1} + y^{k+1}$ for all $k \ge 0$.
\item[(iii)] $f_m \equiv 1 \pmod 2$ if and only if $m + 1$ is a power of $2$. In particular, if $D_s = 1$, then $m = 2^K - 1$ for some $K \ge 1$.
\end{enumerate}
\end{lemma}

\begin{proof}
(i) holds because every $s_e$ is odd. For (ii), both sides are equal to $1$ for $k = 0$ and $k = 1$, and both satisfy the recurrence $g_k = g_{k-1} + (y + y^2)\, g_{k-2}$: the left side by \eqref{eq:rec}, the right side because $y$ and $1 + y$ are the roots of $r^2 + r + (y + y^2)$ over $\F$. For (iii), note that the substitution $x \mapsto y + y^2$ is injective on $\F[x]$. By (ii), $f_m \equiv 1 \pmod 2$ if and only if $(1 + y)^{m+1} = 1 + y^{m+1}$ in $\F[y]$. This holds if $m + 1$ is a power of $2$. If $m + 1 = 2^a b$ with $b > 1$ odd, then $(1 + y)^{m+1} = (1 + y^{2^a})^b$ contains the term $y^{2^a}$, and $2^a < m + 1$. The last assertion follows from (i).
\end{proof}

The last assertion of Lemma~\ref{lem:mod2} is \cite[Theorem~4.1]{Cor21}; the proof given there uses Lucas's theorem. The closed form (ii), which is Binet's formula for the Fibonacci polynomials \cite{HB73} read modulo $2$, will be needed in the next section.

\section{Reduction modulo 4 and proof of the main theorems}\label{sec:mod4}

For $1 \le e \le m - 1$ let $c_e = x\, f_{e-1}\, f_{m-1-e} \in \mathbb{Z}[x]$. Deleting the endpoints of edge $e$ leaves paths with $e - 1$ and $m - 1 - e$ vertices, so $c_e = \sum_{M \ni e} x^{|M|}$, the sum over the matchings of $P_m$ that contain $e$. Write $s_e = 1 - 2t_e$ with $t_e \in \{0, 1\}$. Since $\prod_{e \in M} (1 - 2t_e) \equiv 1 - 2 \sum_{e \in M} t_e \pmod 4$ for every matching $M$, summing over $M$ gives
\begin{equation}\label{eq:lin}
D_s \equiv f_m - 2 \sum_{e=1}^{m-1} t_e\, c_e \pmod 4.
\end{equation}
Clearly $c_e = c_{m-e}$, and for odd $m$ these are the only linear relations among the $c_e$ modulo $2$:

\begin{lemma}\label{lem:indep}
Let $m = 2N + 1$. The reductions of $c_1, \dots, c_N$ modulo $2$ are linearly independent over $\F$.
\end{lemma}

\begin{proof}
Suppose $\sum_{e \in S} c_e \equiv 0 \pmod 2$ for some nonempty $S \subseteq \{1, \dots, N\}$. Since $\F[x]$ is a domain, $\sum_{e \in S} f_{e-1} f_{m-1-e} \equiv 0$. Substituting $x = y + y^2$ and using Lemma~\ref{lem:mod2},
\[
f_{e-1} f_{m-1-e} \mapsto \big((1 + y)^e + y^e\big)\big((1 + y)^{m-e} + y^{m-e}\big) = A + \varphi_e,
\]
where $A = (1 + y)^m + y^m$ has constant term $1$, and $\varphi_e = y^e (1 + y)^{m-e} + y^{m-e} (1 + y)^e$. For $e \le N$ we have $m - e \ge N + 1 > e$, so the lowest-degree term of $\varphi_e$ is $y^e$. Comparing constant terms in $|S|\, A + \sum_{e \in S} \varphi_e = 0$ shows that $|S|$ is even, hence $\sum_{e \in S} \varphi_e = 0$. But the coefficient of $y^{\min S}$ in this sum is $1$, because every other $\varphi_e$ has lowest-degree term of higher degree. This is a contradiction.
\end{proof}

\begin{proposition}\label{prop:mod4}
Let $m = 2N + 1 = 2^K - 1$ with $K \ge 2$, and let $s \in \{\pm 1\}^{m-1}$. Then $D_s \equiv 1 \pmod 4$ if and only if $s$ is folded, that is, $s = (\sigma, \varepsilon, -\varepsilon, \rev\sigma)$ for some $\sigma \in \{\pm 1\}^{N-1}$ and $\varepsilon = \pm 1$. Equivalently, $s_N = -s_{N+1}$, and $s_{m-e} = s_e$ for $1 \le e \le N - 1$.
\end{proposition}

\begin{proof}
By Lemma~\ref{lem:mod2}, $f_m = 1 + 2g$ for some $g \in \mathbb{Z}[x]$. By \eqref{eq:lin}, $D_s \equiv 1 \pmod 4$ is equivalent to
\begin{equation}\label{eq:system}
\sum_{e=1}^{m-1} t_e\, c_e \equiv g \pmod 2.
\end{equation}
The pairs $\{e, m - e\}$ with $1 \le e \le N$ partition $\{1, \dots, 2N\}$, and $c_e = c_{m-e}$. So the left side of \eqref{eq:system} equals $\sum_{e=1}^{N} w_e c_e$ with $w_e = t_e + t_{m-e} \bmod 2$. By Lemma~\ref{lem:indep}, at most one vector $w \in \F^N$ satisfies \eqref{eq:system}.

Folded sequences satisfy the congruence. Indeed, if $s = (\sigma, \varepsilon, -\varepsilon, \rev\sigma)$, then $D_s = D_\sigma^2$ by Lemma~\ref{lem:fold}, and $D_\sigma \equiv f_N \equiv 1 \pmod 2$ by Lemma~\ref{lem:mod2}, because $N + 1 = 2^{K-1}$. Writing $D_\sigma = 1 + 2h$ gives $D_s = 1 + 4h + 4h^2 \equiv 1 \pmod 4$. A folded sequence has $w_e = 0$ for $e < N$, because its two halves are mirror images, and $w_N = 1$, because its two middle signs differ. Hence every $s$ with $D_s \equiv 1 \pmod 4$ has $w = (0, \dots, 0, 1)$, which means that $s$ is folded.
\end{proof}

\begin{remark}\label{rem:mod4}
The congruence places no condition on $\sigma$, so it has $2^N$ solutions. By Lemma~\ref{lem:fold}, a solution satisfies $D_s = 1$ if and only if $D_\sigma = 1$.
\end{remark}

\begin{proof}[Proof of Theorems~\ref{thm:main} and~\ref{thm:main-s}]
By Lemma~\ref{lem:charpoly} it suffices to prove Theorem~\ref{thm:main-s} and the count. Suppose $D_s = 1$. By Lemma~\ref{lem:mod2}, $m = 2^K - 1$ for some $K \ge 1$. We show $s \in F_K$ by induction on $K$. For $K = 1$, $s$ is empty. Let $K \ge 2$. Since $D_s \equiv 1 \pmod 4$, Proposition~\ref{prop:mod4} gives $s = (\sigma, \varepsilon, -\varepsilon, \rev\sigma)$ with $\sigma$ of length $2^{K-1} - 2$. By Lemma~\ref{lem:fold}, $D_\sigma^2 = D_s = 1$, and since $D_\sigma(0) = 1$ this forces $D_\sigma = 1$. By induction $\sigma \in F_{K-1}$, hence $s \in F_K$.

Conversely, if $s \in F_K$, then $D_s = 1$ by Lemma~\ref{lem:fold} and induction; this is \cite[Theorem~4.2]{Cor21}. The count follows from $|F_K| = 2^{K-1}$, because distinct vectors $u$ give distinct matrices $T_m(u)$.
\end{proof}

\begin{example}
For $m = 7$ a folded sequence has $D_s = \big(1 + (s_1 + s_2)x\big)^2$, by Lemma~\ref{lem:fold}. Thus four of the eight folded sequences are nilpotent, namely those with $s_1 = -s_2$. For the others the congruence holds but the equation fails; for instance, $s = (+, +, +, -, +, +)$ has $D_s = 1 + 4x + 4x^2$.
\end{example}

\begin{proof}[Proof of Theorem~\ref{thm:sign}]
Apply Lemma~\ref{lem:charpoly} and Theorem~\ref{thm:main-s} with $s_k = -a_k c_k$, and recall that $F_K$ is closed under negation. The statement on $A_m^b$ is the case $a = (1, \dots, 1)$, $c = b$.
\end{proof}

\begin{corollary}\label{cor:pairs}
If $T_m(u)$ is nilpotent, then each pair $(u_{2j-1}, u_{2j})$ is $(1, i)$ or $(i, 1)$.
\end{corollary}

\begin{proof}
This follows by induction on $K$, since $\sigma$ has even length in $(\sigma, \varepsilon, -\varepsilon, \rev\sigma)$.
\end{proof}

\begin{remark}\label{rem:jordan}
Deleting the first column and the last row of $J_m(a, c)$ leaves a triangular matrix with the non-zero diagonal entries $a_1, \dots, a_{m-1}$. Hence $J_m(a, c)$ has rank at least $m - 1$, and a nilpotent $J_m(a, c)$ consists of a single Jordan block.
\end{remark}

\section{The continued fraction of a nilpotent vector}\label{sec:cf}

All nilpotent matrices $T_m(u)$ of a given dimension are similar (Remark~\ref{rem:jordan}). The $(1,1)$ entry of the resolvent distinguishes them: it is a polynomial in $\lambda^{-1}$ with coefficients in $\{0, \pm 1\}$, whose non-zero coefficients record the choices made in the recursion (Proposition~\ref{prop:cf} and Corollary~\ref{cor:signs}). We compute it by means of a continued fraction.

Expanding at the first vertex instead of the last gives $D_{1,m} = D_{2,m} + s_1 x\, D_{3,m}$. By induction, the terminating continued fraction
\[
\cfrac{1}{1 + \cfrac{s_1 x}{1 + \cfrac{s_2 x}{\ddots + \cfrac{s_{m-1} x}{1}}}}
\]
equals $D_{2,m} / D_s$. For a nilpotent vector the denominator is $1$, so the continued fraction is the polynomial $D_{2,m} = D_{s'}$, where $s'$ is $s$ without its first sign. To determine it we use the following simple identity.

\begin{lemma}\label{lem:cassini}
For $2 \le k \le m$,
\[
D_{1,k}\, D_{2,k-1} - D_{2,k}\, D_{1,k-1} = (-1)^k\, x^{k-1} \prod_{e=1}^{k-1} s_e.
\]
\end{lemma}

\begin{proof}
Let $W_k$ denote the left side. By \eqref{eq:rec}, $W_2 = s_1 x$ and $W_k = -s_{k-1}\, x\, W_{k-1}$ for $k \ge 3$.
\end{proof}

\begin{proposition}\label{prop:cf}
\interlinepenalty=10000 \postdisplaypenalty=10000 % keep the statement on one page
Let $K \ge 2$ and $s \in F_K$, and for $2 \le j \le K$ let $\varepsilon_j$ be the sign $\varepsilon$ chosen in \eqref{eq:rule-s} at the step from $F_{j-1}$ to $F_j$ in the construction of $s$; thus $\varepsilon_j = s_{2^{j-1}-1}$. Then
\begin{equation}\label{eq:lacunary}
D_{s'} = 1 - \varepsilon_2\, x + \sum_{j=3}^{K} \varepsilon_j\, x^{2^{j-1}-1}.
\end{equation}
In particular, $s \mapsto D_{s'}$ is a bijection from $F_K$ onto the $2^{K-1}$ polynomials $1 + \sum_{j=1}^{K-1} \pm x^{2^j - 1}$.
\end{proposition}

\begin{proof}
For $K = 2$ we have $s = (\varepsilon_2, -\varepsilon_2)$ and $D_{s'} = 1 - \varepsilon_2 x$. Let $K \ge 3$ and write $s = (\sigma, \varepsilon_K, -\varepsilon_K, \rev\sigma)$, where $\sigma \in F_{K-1}$ has the signs $\varepsilon_2, \dots, \varepsilon_{K-1}$ and length $N - 1$ with $N = 2^{K-1} - 1$. Let $D_{a,b}$ refer to the path of $s$. The path of $s'$ has the vertices $2, \dots, 2N + 1$ and the middle vertex $c = N + 1$. Its right part has matching polynomial $D_\sigma = D_{1,N} = 1$ by Theorem~\ref{thm:main-s}, and without its first vertex it has matching polynomial $D_{1,N-1}$, by reflection. Classifying the matchings by how $c$ is covered, as in Lemma~\ref{lem:fold}, and applying Lemma~\ref{lem:cassini} with $k = N$, which is odd, gives
\[
D_{s'} = D_{2,N} + \varepsilon_K x \big(D_{1,N}\, D_{2,N-1} - D_{2,N}\, D_{1,N-1}\big) = D_{\sigma'} - \varepsilon_K x^N \prod_{e=1}^{N-1} \sigma_e,
\]
where $\sigma'$ is $\sigma$ without its first sign. Since $\sigma$ is folded, say $\sigma = (\tau, \varepsilon', -\varepsilon', \rev\tau)$, we have $\prod_e \sigma_e = -\big(\prod_e \tau_e\big)^2 = -1$, so $D_{s'} = D_{\sigma'} + \varepsilon_K x^N$, and \eqref{eq:lacunary} follows by induction on $K$. Since $s$ is determined by $\varepsilon_2, \dots, \varepsilon_K$, the map $s \mapsto D_{s'}$ is a bijection.
\end{proof}

\begin{corollary}\label{cor:signs}
Let $T = T_m(u)$ be nilpotent, with $m = 2^K - 1 \ge 3$, and let $\varepsilon_j$ be the signs of Proposition~\ref{prop:cf}. Then
\begin{equation}\label{eq:resolvent}
\big((\lambda I - T)^{-1}\big)_{11} = \lambda^{-1} - \varepsilon_2\, \lambda^{-3} + \sum_{j=3}^{K} \varepsilon_j\, \lambda^{-(2^j-1)},
\end{equation}
so $(T^k)_{11} = 0$ for every $k \ge 0$ that is not of the form $2^j - 2$ with $1 \le j \le K$.
\end{corollary}

\begin{proof}
By Cramer's rule and Lemma~\ref{lem:charpoly}, applied to $T$ and to $T$ without its first row and column, $\big((\lambda I - T)^{-1}\big)_{11} = \lambda^{-1} D_{s'}(\lambda^{-2})$, and \eqref{eq:lacunary} gives \eqref{eq:resolvent}. The last claim follows by comparing coefficients in $(\lambda I - T)^{-1} = \sum_{k=0}^{m-1} T^k \lambda^{-k-1}$.
\end{proof}

\ifwithdoubling
\section{A doubling rule}\label{sec:doubling}

By Corollary~\ref{cor:pairs}, a nilpotent sign sequence has the form $(\theta_1, -\theta_1, \dots, \theta_N, -\theta_N)$. The next identity leads to a second recursive description of $F_K$, in which new signs are interleaved with old ones.

\begin{proposition}\label{prop:double}
Let $N \ge 1$ and $\theta_1, \dots, \theta_N \in \{\pm 1\}$. Then
\[
D_{(\theta_1, -\theta_1, \theta_2, -\theta_2, \dots, \theta_N, -\theta_N)}(x) = D_{(\theta_1\theta_2,\, \theta_2\theta_3,\, \dots,\, \theta_{N-1}\theta_N)}(x^2).
\]
\end{proposition}

\begin{proof}
Let $r = (\theta_1, -\theta_1, \dots, \theta_N, -\theta_N)$ and $\tau = (\theta_1\theta_2, \dots, \theta_{N-1}\theta_N)$, and put $Q_j = D_{1,2j+1}$, where $D_{a,b}$ refers to the path of $r$. Then $Q_0 = Q_1 = 1$, and applying \eqref{eq:rec} three times gives, for $2 \le j \le N$,
\[
Q_j = \big(1 + (r_{2j-1} + r_{2j})\, x\big) Q_{j-1} - r_{2j-2}\, r_{2j-1}\, x^2\, Q_{j-2} = Q_{j-1} + \tau_{j-1}\, x^2\, Q_{j-2},
\]
because $r_{2j-1} + r_{2j} = 0$ and $r_{2j-2}\, r_{2j-1} = -\theta_{j-1}\theta_j$. This is the recurrence \eqref{eq:rec} for the path with the signs $\tau$, with $x$ replaced by $x^2$. Hence $D_r(x) = Q_N(x) = D_\tau(x^2)$.
\end{proof}

\begin{corollary}\label{cor:double}
Let $K \ge 2$ and $N = 2^{K-1} - 1$. For $\eta = \pm 1$ and $\tau \in \{\pm 1\}^{N-1}$ put
\[
\Phi_\eta(\tau) = (\theta_1, -\theta_1, \dots, \theta_N, -\theta_N), \qquad \theta_1 = \eta, \quad \theta_{j+1} = \tau_j\, \theta_j.
\]
Then $(\eta, \tau) \mapsto \Phi_\eta(\tau)$ is a bijection from $\{\pm 1\} \times F_{K-1}$ onto $F_K$.
\end{corollary}

\begin{proof}
Since $\theta_j\theta_{j+1} = \tau_j$, Proposition~\ref{prop:double} gives $D_{\Phi_\eta(\tau)}(x) = D_\tau(x^2)$, so by Theorem~\ref{thm:main-s}, $\Phi_\eta(\tau) \in F_K$ if and only if $\tau \in F_{K-1}$. By Corollary~\ref{cor:pairs}, every $s \in F_K$ is of the form $\Phi_\eta(\tau)$, where $\eta = s_1$ and $\tau_j = s_{2j-1}\, s_{2j+1}$ are determined by $s$.
\end{proof}

\begin{remark}\label{rem:ctilde}
The sequence \eqref{eq:ctilde} is a fixed point of this rule in the following sense. Put $b^{(K)} = (\tc_2, \dots, \tc_{2^K - 1})$ and $\theta_j = -\tc_{2j}$. Then \eqref{eq:ctilde} gives $-\tc_{2j+1} = -\theta_j$, $\theta_1 = -1$ and $\theta_j \theta_{j+1} = -\tc_{j+1}$, that is, $-b^{(K)} = \Phi_{-1}\big(-b^{(K-1)}\big)$ for $K \ge 2$. By Proposition~\ref{prop:double}, $D_{-b^{(K)}}(x) = D_{-b^{(K-1)}}(x^2) = \dots = 1$, so $A_m^{b^{(K)}}$ is nilpotent for $m = 2^K - 1$. This gives another proof of \cite[Lemma~30]{CWD12}.
\end{remark}

\begin{remark}
In terms of matrices, Proposition~\ref{prop:double} reflects the fact that, with $r$ and $\tau$ as in its proof, the square of $A_{2N+1}^{-r}$ maps the span of the even-numbered basis vectors into itself and acts on it as $A_N^{-\tau}$. This is the finite counterpart of \cite[Lemma~5]{CWD12}, from which Chandler-Wilde and Davies deduce that the spectrum $S$ of the Feinberg--Zee hopping matrix is closed under taking square roots: if $\lambda \in S$, then $\pm\sqrt{\lambda} \in S$. The maps $\Phi_{\pm 1}$ correspond to their maps $\Gamma_\pm$.
\end{remark}
\fi

\section{Formal verification and computations}\label{sec:lean}

We have formalised the proof of Theorem~\ref{thm:main} in Lean~4 \cite{dMU21}, version 4.36.0, using the Mathlib library \cite{mathlib}. The formalisation follows Sections~\ref{sec:matching} and~\ref{sec:mod4} and has about $1{,}150$ lines. The matrix $T_m(u)$, with rows and columns indexed from $0$, is defined by
\begin{Verbatim}
def T (m : ℕ) (u : ℕ → ℂ) : Matrix (Fin m) (Fin m) ℂ :=
  fun i j => if (j : ℕ) = i + 1 then u i else if (i : ℕ) = j + 1 then -u j else 0
\end{Verbatim}
\noindent and the vectors of Theorem~\ref{thm:main} by the inductive predicate
\begin{Verbatim}
inductive CorlessVec : ℕ → List ℂ → Prop
  | base : CorlessVec 0 []
  | step₁ {k : ℕ} {u : List ℂ} :
      CorlessVec k u → CorlessVec (k + 1) (u ++ [1, Complex.I] ++ u.reverse)
  | step₂ {k : ℕ} {u : List ℂ} :
      CorlessVec k u → CorlessVec (k + 1) (u ++ [Complex.I, 1] ++ u.reverse)
\end{Verbatim}
\noindent in which the index $k$ corresponds to $K - 1$. The formal statements of Theorem~\ref{thm:main} are
\begin{Verbatim}
theorem corless_conjecture {u : List ℂ} (hu : ∀ x ∈ u, x = 1 ∨ x = Complex.I) :
    IsNilpotent (T (u.length + 1) (fun i => u.getD i 1)) ↔ ∃ k, CorlessVec k u
\end{Verbatim}
\begin{Verbatim}
theorem card_nilpotent (k : ℕ) :
    {u : List ℂ | u.length = 2 ^ (k + 1) - 2 ∧ (∀ x ∈ u, x = 1 ∨ x = Complex.I) ∧
      IsNilpotent (T (u.length + 1) (fun i => u.getD i 1))}.ncard = 2 ^ k
\end{Verbatim}
\noindent The first of them includes the restriction on the dimension, since every vector satisfying \texttt{CorlessVec k} has length $2^{k+1} - 2$. Both theorems depend only on the three standard axioms of Lean (\texttt{propext}, \texttt{Classical.choice} and \texttt{Quot.sound}). Theorem~\ref{thm:sign} and \furthersections{} have not been formalised.

As an independent check, we have verified the case $m = 63$ of Theorem~\ref{thm:main} by an exhaustive search. Every nilpotent sign sequence satisfies $D_s \equiv 1 \pmod 4$ and therefore solves the linear system \eqref{eq:system}, which has rank $31$ for $m = 63$. This leaves $2^{31}$ of the $2^{62}$ sign sequences. For each of them we evaluated $D_s$ at two points modulo the prime $2^{61} - 1$, and we recomputed $D_s$ exactly whenever both values were equal to $1$. Exactly $32$ sequences have $D_s = 1$, namely the members of $F_6$. The search takes about $20$ seconds on $14$ processor cores. We have also verified Theorem~\ref{thm:sign} for $m \le 9$, and the results of \furthersections{} for $K \le 8$, by exact computation. The formalisation and the programs are available at \url{https://github.com/marcelbtec/corless-nilpotent} and are archived at \cite{Bla26}.

\section{Concluding remarks}

Theorem~\ref{thm:sign} with $a = u$ and $c = -u$ covers the population $\{1, -1, i, -i\}$: since each sign sequence comes from $2^{m-1}$ vectors $u$, there are $2^{K-1} \cdot 2^{m-1}$ nilpotent matrices of dimension $m = 2^K - 1$. If $0$ is admitted as an entry as well, the zero entries divide $u$ into (possibly empty) blocks $v$, the matrix $T_m(u)$ is the direct sum of the corresponding matrices $T_{|v|+1}(v)$, where $|v|$ is the length of $v$, and it is nilpotent if and only if each of them is.

By Lemma~\ref{lem:charpoly}, the multiplicity of the eigenvalue $0$ of $T_m(u)$ is $m - 2 \deg D_s$. Corless observed that for $m = 15$ the multiplicities $3, 5, 7, 11, 13$ and $15$ occur, but $9$ does not \cite[Section~5]{Cor21}; in other words, $\deg D_s \ne 3$ for all $s \in \{\pm 1\}^{14}$. It would be interesting to determine the degrees that $D_s$ can have for general $m$. Theorem~\ref{thm:main-s} settles the case $\deg D_s = 0$, and conditions of the form $\deg D_s \le d$ might be accessible by the $2$-adic method used here.

\section*{Acknowledgements}
This work grew out of the paper \cite{Cor21} of R.~M.~Corless, which poses the problem and contains several ingredients of its solution.

\section*{Collaboration with an AI system}
This paper is the result of a collaboration between the author and the AI system Claude (\mbox{Anthropic}; principally Opus~5.5 and Fable~5.2). The system  wrote the Lean formalisation and the programs of Section~\ref{sec:lean}; the author contributed the sketch for the proof of Theorem~\ref{thm:main}, Theorem~\ref{thm:sign} and the results of \furthersections{} and wrote the second version of the paper, which the system checked and both revised. The author has checked all arguments and takes full responsibility for the paper.


\begin{thebibliography}{99}

\bibitem{Bla26}
M.~Blattner, \emph{Corless's nilpotent conjecture: proof, Lean formalisation and computations}, version 1.1.0, Zenodo, 2026. \url{https://doi.org/10.5281/zenodo.24317787}.

\bibitem{CCL11}
S.~N.~Chandler-Wilde, R.~Chonchaiya and M.~Lindner, \emph{Eigenvalue problem meets Sierpinski triangle: computing the spectrum of a non-self-adjoint random operator}, Operators and Matrices \textbf{5} (2011), no.~4, 633--648. \url{https://doi.org/10.7153/oam-05-46}.

\bibitem{CCL13}
S.~N.~Chandler-Wilde, R.~Chonchaiya and M.~Lindner, \emph{On the spectra and pseudospectra of a class of non-self-adjoint random matrices and operators}, Operators and Matrices \textbf{7} (2013), no.~4, 739--775. \url{https://doi.org/10.7153/oam-07-43}.

\bibitem{CWD12}
S.~N.~Chandler-Wilde and E.~B.~Davies, \emph{Spectrum of a Feinberg--Zee random hopping matrix}, Journal of Spectral Theory \textbf{2} (2012), no.~2, 147--179. \url{https://doi.org/10.4171/JST/25}; arXiv:1110.0792.

\bibitem{Cor21}
R.~M.~Corless, \emph{Skew-symmetric tridiagonal Bohemian matrices}, Maple Transactions \textbf{1} (2021), no.~2, Article 14360. \url{https://doi.org/10.5206/mt.v1i2.14360}.

\bibitem{CLPS22}
R.~M.~Corless, G.~Labahn, D.~Piponi and L.~Rafiee Sevyeri, \emph{Bohemian matrix geometry}, in: Proceedings of the 2022 International Symposium on Symbolic and Algebraic Computation (ISSAC 2022), ACM, 2022, 361--370. \url{https://doi.org/10.1145/3476446.3536177}; arXiv:2202.07769.

\bibitem{dMU21}
L.~de~Moura and S.~Ullrich, \emph{The Lean 4 theorem prover and programming language}, in: Automated Deduction -- CADE 28, Lecture Notes in Computer Science \textbf{12699}, Springer, 2021, 625--635. \url{https://doi.org/10.1007/978-3-030-79876-5_37}.

\bibitem{FZ99}
J.~Feinberg and A.~Zee, \emph{Non-Hermitian localization and delocalization}, Physical Review E \textbf{59} (1999), no.~6, 6433--6443. \url{https://doi.org/10.1103/PhysRevE.59.6433}.

\bibitem{Hag15a}
R.~Hagger, \emph{Symmetries of the Feinberg--Zee random hopping matrix}, Random Matrices: Theory and Applications \textbf{4} (2015), no.~4, 1550016. \url{https://doi.org/10.1142/S2010326315500161}.

\bibitem{Hag15b}
R.~Hagger, \emph{The eigenvalues of tridiagonal sign matrices are dense in the spectra of periodic tridiagonal sign operators}, Journal of Functional Analysis \textbf{269} (2015), no.~5, 1563--1570. \url{https://doi.org/10.1016/j.jfa.2015.01.019}; arXiv:1412.1724.

\bibitem{HB73}
V.~E.~Hoggatt, Jr. and M.~Bicknell, \emph{Roots of Fibonacci polynomials}, Fibonacci Quarterly \textbf{11} (1973), no.~3, 271--274. \url{https://doi.org/10.1080/00150517.1973.12430825}.

\bibitem{mathlib}
The mathlib Community, \emph{The Lean mathematical library}, in: Proceedings of the 9th ACM SIGPLAN International Conference on Certified Programs and Proofs (CPP 2020), ACM, 2020, 367--381. \url{https://doi.org/10.1145/3372885.3373824}.

\end{thebibliography}
\end{document}